\documentclass[12pt,leqno]{article}
\usepackage{amsmath, amsthm, amsfonts, amssymb, color}
\usepackage{mathrsfs}
\theoremstyle{definition}

\newcommand{\scr}[1]{\mathscr #1}
\definecolor{wco}{rgb}{0.5,0.2,0.3}

\numberwithin{equation}{section} \theoremstyle{remark}

\newcommand{\ua}{\uparrow}

\title{
{\bf Log-Harnack Inequality for Mild Solutions of SPDEs with Strongly Multiplicative Noise}
\footnote{Supported in part by NNSFC(11131003), SRFDP, the Laboratory of Mathematical and  Complex Systems and the Fundamental Research Funds for the Central Universities.}
}
\author{
{\bf Feng-Yu Wang$^{1),3)}$, Tusheng Zhang$^{2)}$ }\\ \\
\footnotesize{$^1)$ School of Mathematical Sciences, Beijing Normal University, Beijing 100875, China}\\
 \footnotesize{$^2)$ School of Mathematics, University of Manchester, Oxford Road, Manchester M13 9PL, UK} \\
  \footnotesize{$^3)$ Department of Mathematics, Swansea University, Singleton Park, SA2 8PP, UK}\\
\footnotesize{Email: \tttext{wangfy@bnu.edu.cn}; \tttext{tusheng.zhang@manchester.ac.uk}} }

\begin{document}
\def\tttext#1{{\normalfont\ttfamily#1}}
\def\R{\mathbb R}  \def\ff{\frac} \def\ss{\sqrt} \def\B{\mathbf B}
\def\N{\mathbb N} \def\kk{\kappa} \def\m{{\bf m}}
\def\dd{\delta} \def\DD{\Delta} \def\vv{\varepsilon} \def\rr{\rho}
\def\<{\langle} \def\>{\rangle} \def\GG{\Gamma} \def\gg{\gamma}
  \def\nn{\nabla} \def\pp{\partial} \def\EE{\scr E}
\def\d{\text{\rm{d}}} \def\bb{\beta} \def\aa{\alpha} \def\D{\scr D}
  \def\si{\sigma} \def\ess{\text{\rm{ess}}}
\def\beg{\begin} \def\beq{\begin{equation}}  \def\F{\scr F}
\def\Ric{\text{\rm{Ric}}} \def\Hess{\text{\rm{Hess}}}
\def\e{\text{\rm{e}}} \def\ua{\underline a} \def\OO{\Omega}  \def\oo{\omega}
 \def\tt{\tilde} \def\Ric{\text{\rm{Ric}}}
\def\cut{\text{\rm{cut}}} \def\P{\mathbb P}
\def\C{\scr C}     \def\E{\mathbb E}
\def\Z{\mathbb Z} \def\II{\mathbb I}
  \def\Q{\mathbb Q}  \def\LL{\Lambda}\def\L{\scr L}
  \def\B{\scr B}    \def\ll{\lambda}
\def\vp{\varphi}\def\H{\mathbb H}\def\ee{\mathbf e}

\maketitle
\begin{abstract} Due to technical reasons, existing results concerning   Harnack type inequalities for SPDEs with multiplicative noise apply only to  the case where the coefficient in the noise term  is an Hilbert-Schmidt perturbation of a fixed bounded operator. In this paper we investigate a class of semi-linear SPDEs with strongly multiplicative noise whose coefficient is  even allowed to be unbounded which is thus no way to be Hilbert-Schmidt. Gradient estimates, log-Harnack inequality and applications are derived. Applications to stochastic reaction-diffusion equations driven by space-time white noise are presented. \end{abstract} \noindent

 \vskip 2cm

\section{Introduction}

Let $(\H,\<\cdot,\cdot\>,|\cdot|)$ be a separable Hilbert space.   Let $\L(\H)$ be the set of all densely defined linear operators on $\H$. We will use $\|\cdot\|$ and $\|\cdot\|_{HS}$ to denote the operator norm and the Hilbert-Schmidt norm for linear operators on $\H$ respectively. Let
$$b: \H\to \H,\ \ \si: \H\to \L(\H)$$ be two given measurable maps.

 Consider the following SPDE on $\H$:
\beq\label{E} \d X_t= \big\{AX_t+ b(X_t)\big\}\d t +\si(X_t)\d W_t,\end{equation}
where $W_t, t\geq 0$ is a cylindrical Brownian motion on $H$ admitting the representation:
\begin{equation}\label{Cylin}
W_t=\sum_{i=1}^{\infty}\beta_i(t)e_i.
\end{equation}
Here $\beta_i(t), i\geq 1$ is a sequence of independent real-valued Brownian motions on a complete filtered probability space $(\OO, (\F_t)_{t\ge 0},\P).$

Recall that an $\H$-valued adapted process $(X_t)_{t\ge 0}$ is called a mild solution to (\ref{E}) if
\beq\label{Sol}  \int_0^t \E(|T_{t-s}b(X_s)|+\|T_{t-s} \si(X_s)\|_{HS}^2)\d s<\infty,\ \ t\ge 0\end{equation}  and almost surely
$$X_t = T_t X_0 + \int_0^t T_{t-s} b(X_s)\d s +\int_0^tT_{t-s}\si(X_s)\d W_s,\ \ t\ge 0.$$

To ensure the existence and uniqueness of the mild solution, and to derive regularity estimates of the associated semigroup, we shall make use of the following conditions.
\beg{enumerate} \item[{\bf (A1)}] There exists a positive function $K_b\in C((0,\infty))$ such that
$$\phi_b(t):= \int_0^t K_b(s)\d s<\infty,\  |T_t( b(x)-b(y))|^2\le K_b(t)|x-y|^2,\ \  t>0, x,y\in\H.$$
\item[{\bf (A2)}] $|\si v|^2\ge \ll(\si)|v|^2$ holds  for some constant  $\ll(\si) >0$ and all $v\in\H.$
\item[{\bf (A3)}] There exists $x\in\H$ such that for any $s>0,$ $T_s\si(x)$ extends to an unique Hilbert-Schmidt operator which is again denote by $T_s\si(x)$ such that
$\int_0^t\|T_s\si(x)\|_{HS}^2\d s<\infty, t>0$; and   there exists a positive measurable function $K_\si\in C((0,\infty))$ such that
$$ \phi_\si(t):=\int_0^tK_\si(s)\d s<\infty,\   \|T_t(\si(x)-\si(y)\|_{HS}^2 \le K_\si(t)|x-y|^2,\ \  t>0, x,y\in\H.$$
\item[{\bf (A4)}] The operator $A$ admits a complete orthonormal system of eigenvectors, that is, there exists an orthonormal basis $\{e_n, n\geq 1\}$ of $\H$ such that $-Ae_n=\lambda_ne_n, n\geq 1$, where $ \lambda_n\geq 0, n\geq 0$ are the corresponding eigenvalues.
\end{enumerate}

Replacing (A.3) and (B.1) in the proof of \cite[Theorem A.1]{PZ} by the current  {\bf (A1)} and {\bf (A3)}, we see    that  for any $X_0\in L^2(\OO,\F_0,\P)$  the equation has a unique mild solution
and  $\E |X_t|^2$ is locally bounded in $t$.
 Let
$$P_t f(x)= \E f(X_t^x),\ \ x\in\H, f\in\B_b(\H), t\ge 0,$$ where $X_t^x$ denotes the unique mild solution to (\ref{E}) with $X_0=x.$ Under {\bf (A1)}-{\bf (A3)}   the proof of  \cite[Theorem 1.2]{PZ} implies  that  $P_tf$ is Lipschitz continuous for any $t>0$ and $f\in\B_b(\H);$ consequently, $P_t$ is strong Feller.

\

In this paper we aim to investigate Harnack type inequalities for $P_t$, which implies not only the strong Feller property but also some concrete estimates on the heat kernel.   As the process is infinite-dimensional, the Harnack inequality we shall establish will be dimension-free. The following type dimension-free Harnack inequality with a power $\aa>1$
$$|P_tf|^\aa(x)\le P_t |f|^\aa(y) \e^{C(t)\rr(x,y)^2}$$ was first found in \cite{W97} for diffusion semigroups on Riemannian manifolds, where $\rr$ is the Riemannian distance. Because of the new coupling argument introduced in \cite{ATW06}, this inequality has been established for a large class of SDEs and SPDEs (see \cite{GW,L,LW, O, W11,WY11,ZS,ZT}   and references therein). When this type of inequality is invalid, the following weaker version, known as log-Harnack inequality, was investigated as a substitution (see \cite{RW10, W10, WWX10, X}):
$$P_t\log f(x)\le \log P_t f(y) + C(t)|x-y|^2,\ \ t>0, x,y\in \H,\  \ f>0, f\in \B_b(\H).$$
However,   when SPDEs with  multiplicative noise is considered, existing results on Harnack type inequalities work only for the case that the coefficient  in the noise term is   an Hilbert-Schmidt perturbation of a constant operator; i.e. the conditions imply that $\|\si(x)-\si(y)\|_{HS}<\infty$ for $x,y\in\H$. Although this assumption comes out naturally by applying the It\^o's formula to the distance of the two martingale processes of the coupling, it however excludes many important models;
for instance  $A$ being the Dirichlet Laplacian on $[0,1]$ and $\si(x) =(\phi\circ x) {\rm Id}$ for a Lipschitz function $\phi$ on $\R$ as studied in \cite{DP, ZT} where a reflection is also considered (see Section 4 for details).

To get ride of the condition on $\|\si(x)-\si(y)\|_{HS}$, we will not make use of the coupling method, but follow the line of \cite{RW10} by establishing the gradient estimate of type $|\nn P_t f|^2\le C(t)P_t|\nn f|^2$, which, along with a finite-dimensional approximation argument,  will enable us to derive the log-Harnack inequality. Here, for any function $f$ on $\H$ and $x\in \H$, we let $$|\nn f|(x)=\limsup_{y\to x} \ff{|f(y)-f(x)|}{|x-y|}.$$  Moreover, it is easy to see that  {\bf (A1)} and {\bf (A3)} imply  $$t_0:=\sup\Big\{t>0:  \phi_b(t)+ \phi_\si(t)\le\ff 1 6\Big\} >0.$$ The following is the main result of the paper, where when $t_0=\infty,$ we set
$$t_0(6^{\ff t{t_0}}-1)= t_0(1-6^{-\ff t{t_0}})=\lim_{r\to\infty} r(6^{\ff t{r}}-1)=\lim_{r\to\infty} r(1-6^{-\ff t{r}} )=  t\log 6.$$

\beg{thm} \label{T1.1} Assume   {\bf (A1)}, {\bf (A3)} and  {\bf (A4)}. \beg{enumerate} \item[$(1)$] For any $f\in C_b^1(\H)$,  $$|\nn P_t f|^2\le 6^{1+\ff t {t_0}} P_t|\nn f|^2,\ \  t\ge 0.$$    \item[$(2)$] If  {\bf (A2)} holds, then for any strictly positive $f\in \B_b(\H)$,
$$ P_t\log f(y)\le \log P_t f(x) +\ff{3\log 6}{\ll(\si)t_0( 1-6^{-\ff t{t_0}} )}|x-y|^2,\ \ x,y\in\H, t>0.$$
\item[$(3)$]  If    {\bf (A2)} holds, then for any $f\in\B_b(\H)$,
$$|\nn P_tf|^2\le \ff{3\log 6}{t_0\ll(\si)(1- 6^{-\ff t{t_0}} )}\big\{P_t f^2-(P_tf)^2\big\},\ \ t>0.$$
\item[$(4)$] If       $|\si v|^2\le \bar\ll(\si)|v|^2$   holds for some constant $\bar\ll(\si)>0$ and all $v\in\H$, then
$$P_t f^2 -(P_t f)^2 \le \ff{12\bar\ll(\si) t_0(6^{\ff t{t_0}}-1)}{\log 6} P_t|\nn f|^2,\ \ f\in  C_b^1(\H), t\ge 0.$$ \end{enumerate}
 \end{thm}

As application of Theorem \ref{T1.1}, (3) implies that $P_t$ sends bounded measurable functions to Lipschitz continuous functions and is thus strong Feller; (4) provides a Poincar\'e inequality for $P_t$; and the log-Harnack inequality in (3) implies the following assertions on the  quasi-invariant measure and heat kernel estimates (see   Corollary 1.2 in \cite{WY11,RW10,WWX10}).
Recall that a $\si$-finite measure $\mu$ is called quasi-invariant for $\mu$ if $\mu P_t$ is absolutely continuous with respect to $\mu$.

\beg{cor}\label{C1.2} Assume {\bf (A1)}-{\bf (A4)}. Let $\mu$ be a quasi-invariant measure of $(P_t)_{t>0}$. Then
\beg{enumerate} \item[$(1)$] $P_t$ has a density $p_t$ with respect to $\mu$ and
$$\int_\H p_t(x,z)\log \ff{p_t(x,z)}{p_t(y,z)}\, \mu(\d z)\le \ff{3\log 6}{\ll(\si)t_0( 1-6^{-\ff t{t_0}} )}|x-y|^2,\ \ x,y\in\H, t>0.$$
\item[$(2)$] For any $x,y\in \H$ and $t>0$,
$$\int_\H p_t(x,z)p_t(y,z)\,\mu(\d z)\ge \exp\bigg[-\ff{3\log 6}{\ll(\si)t_0( 1-6^{-\ff t{t_0}} )}|x-y|^2\bigg].  $$
\item[$(3)$] If $\mu$ is an invariant probability measure of $P_t$, then $\mu$ has full support and it is the unique invariant probability measure. Moreover, letting $P_t^*$ be the adjoint operator of $P_t$ in $L^2(\mu)$ and let $W_2$ be the quadratic Wasserstein distance with respect to $|\cdot|$,   the following entropy-cost inequality holds:
    $$\mu((P_t^*f)\log P_t^*f)\le \ff{3\log 6}{\ll(\si)t_0( 1-6^{-\ff t{t_0}} )} W_2(f\mu,\mu)^2,\   t>0, f\ge 0, \mu(f)=1.$$\end{enumerate} \end{cor}

To apply Corollary \ref{C1.2}, in particular the third assertion, we need to verify the existence of the invariant probability measure of $P_t$. This can be done by using e.g. \cite[Theorem 6.1.2]{DP2} (see the proof of Theorem \ref{T4.1}(3) below).

The rest of the paper is organized as follows. In Section 2, we establish the finite dimensional approximations to the mild solutions of the SPDEs.  Section 3 is devoted to the proofs of Theorem 1.1 and Corollary 1.2. In Section 4, we apply our results to stochastic reaction-diffusion equations driven by space-time white noise.
\section{ Finite dimensional approximations}
In this section, we will prove a finite dimensional approximation result for the mild solution of equation (\ref{E}) which will be used later. Let $\{e_n, n\geq 1\}$ be the eigenbasis of the operator $A$. Set $\H_n={\rm span}\{e_1, e_2, ..., e_n\}$. Denote by $\scr P_n$ the projection operator from $\H$ into $\H_n$. Note that $\scr P_n$ commutes with the semigroup $T_t, t\geq 0$.  Define for $x=\sum_{i=1}^n\<x,e_i\>e_i\in \H_n$,
$$A_nx=-\sum_{i=1}^n\lambda_i\<x,e_i\>e_i.$$
Then $A_n$ is a bounded linear operator on $\H_n$. Introduce
\begin{equation}\label{2.1}
b_n(x)=\scr P_nb(x), \quad \sigma_n(x) y =\scr P_n(\sigma(x) y ), \quad x,y\in \H.
\end{equation}
Consider the following system of stochastic differential equations in $\H_n$:
\begin{equation}\label{2.2}
\left\{\begin{array}{ll}
\d X_t^n& = A_nX_t^n\d t+ b_n(X_t^n)\d t +\si_n(X_t^n)\d W_t^n \\
X_0^n&=\scr P_nX_0,\end{array}\right.
\end{equation}
where $W_t^n=\sum_{i=1}^n\beta_i(t)e_i$. It is well known that under {\bf (A1)} and {\bf  (A3)} the above equation admits a unique strong solution.

\beg{thm}\label{TF}
Let $X_t^n, X_t$ be the $($mild$)$ solutions to equation $(\ref{2.2})$ and $(\ref{E})$. Assume {\bf (A1)}, {\bf (A3)} and {\bf (A4)}. If $\E |X_0|^2<\infty$ then
\beg{equation}\label{2.3}
\lim_{n\to \infty}\E|X_t^n-X_t|^2=0,\ \ t\ge 0.
\end{equation}
\end{thm}
\beg{proof} Fix an arbitrary positive constant $T>0$. We will prove (\ref{2.3}) for $t\leq T$. Let $T^n_t$ denote the semigroup generated by $A_n$. The following representation holds:
$$ T_t^nx=\sum_{i=1}^n\e^{-\lambda_it}\<x,e_i\>e_i, \quad x\in \H_n.$$
In a mild form, we have
\beg{equation}\label{2.4}
X_t^n = T_t^nX_0^n + \int_0^t T_{t-s}^n b_n(X_s^n)\d s +\int_0^tT_{t-s}^n\si_n(X_s^n)\d W_s^n.
\end{equation}
Subtracting $X$ from $X^n$ and taking expectation we get
\beq\label{2.5}\beg{split}
\E|X_t^n-X_t|^2&\leq  3\E |T_t^nX_0^n-T_tX_0|^2  \\
&+3\E \int_0^t|T_{t-s}^n b_n(X_s^n)-T_{t-s}b(X_s)|^2\d s  \\
&+3\E \int_0^t\|T_{t-s}^n\si_n(X_s^n)-T_{t-s}\si(X_s)\|_{HS}^2\d s.\end{split}
\end{equation}
Now, since $T_n^n X_0^n=\scr P_n T_tX_0$ and $T_{t-s}^n b_n= \scr P_n T_{t-s} b,$ we have
\beg{equation}\label{2.6}
\E |T_t^nX_0^n-T_t X_0|^2 =\E \sum_{k=n+1}^{\infty}\<X_0, e_k\>^2\e^{-2\lambda_kt} \leq \E \sum_{k=n+1}^{\infty}\<X_0, e_k\>^2,
\end{equation}
and
\beq\label{2.7}\beg{split}
& \E \int_0^t|T_{t-s}^n b_n(X_s^n)-T_{t-s}b(X_s)|^2\d s  \\
&\leq  2\E[\int_0^t|T_{t-s}^n b_n(X_s^n)-T_{t-s}^nb_n(X_s)|^2\d s +
2\E \int_0^t|T_{t-s}^nb_n(X_s)-T_{t-s}b(X_s)|^2\d s \\
&\leq 2\E \int_0^t|T_{t-s}b(X_s^n)-T_{t-s}b(X_s)|^2\d s +
2\E \int_0^T\sum_{k=n+1}^{\infty}\e^{-2\lambda_k(t-s)}\<b(X_s), e_k\>^2\d s \\
&\leq 2\E \int_0^tK_b(t-s)|X_s^n-X_s|^2\d s +2 \E \int_0^T\sum_{k=n+1}^{\infty}\e^{-2\lambda_k(t-s)}\<b(X_s), e_k\>^2\d s.
\end{split}\end{equation}

To get an upper  bound for the last term in (\ref{2.5}), we observe that
\beq\label{2.8} \beg{split}
&  \E \int_0^t\|T_{t-s}^n \si_n(X_s^n)-T_{t-s}^n\si_n(X_s)\|_{HS}^2\d s =\E\int_0^T\big\|T_{t-s}^n \big(\si(X_s^n)-\si(X_s)\big)\big\|_{HS}^2\d s\\
&\leq   \E \int_0^t\|T_{t-s}\si(X_s^n)-T_{t-s}\si(X_s)\|_{HS}^2\,\d s
 \leq  \E \int_0^tK_{\si}(t-s)|X_s^n-X_s|^2\d s,
\end{split}\end{equation}
where condition {\bf (A3)} was used. Moreover,
\beq \label{2.9}\beg{split}
&  \E \int_0^t\|T_{t-s}^n \si_n(X_s)-T_{t-s}\si(X_s)\|_{HS}^2\,\d s
 = \E \int_0^t\big\|(\scr P_n-\text{Id}) T_{t-s}\si (X_s)\big\|_{HS}^2 \d s \\
&= \E \int_0^t\sum_{m=1}^{\infty} \sum_{k=n+1}^{\infty} \e^{-2\lambda_k(t-s)}\<\si(X_s)e_m, e_k\>^2\d s.
\end{split}\end{equation}
It follows from (\ref{2.8}) and (\ref{2.9}) that
\beq\label{2.10}\beg{split}
&  \E \int_0^t\|T_{t-s}^n\si_n(X_s^n)-T_{t-s}\si(X_s)\|_{HS}^2\,\d s \\
&\leq   2\E \int_0^tK_{\si}(t-s)|X_s^n-X_s|^2\d s
  +2\E \int_0^t\sum_{m=1}^{\infty} \sum_{k=n+1}^{\infty} \e^{-2\lambda_k(t-s)}\<\si(X_s)e_m, e_k\>^2\d s.\end{split}
\end{equation}
Putting (\ref{2.5}), (\ref{2.6}), (\ref{2.7}), (\ref{2.10}) together we arrive at
\beq\label{2.11}\beg{split}
& \E[|X_t^n-X_t|^2] \\
&\leq  C\big\{a_n+c_n(t), d_n(t)\big\}+C\int_0^t(K_{\si}(t-s)+K_b(t-s))\E |X_s^n-X_s|^2 \d s \end{split}
\end{equation} for some constant $C>0$ and
\beg{equation*}\beg{split}
&a_n:=\E \sum_{k=n+1}^{\infty}\<X_0, e_k\>^2,\\
 &c_n(t):=\E \int_0^t\sum_{k=n+1}^{\infty}\e^{-2\lambda_k(t-s)}\<b(X_s), e_k\>^2\d s,\\
&d_n(t):=\E \int_0^t\sum_{m=1}^{\infty} \sum_{k=n+1}^{\infty} \e^{-2\lambda_k(t-s)}\<\si(X_s)e_m, e_k\>^2\d s.\end{split}\end{equation*}
By (\ref{Sol})
and the dominated convergence theorem we see that   $a_n+c_n(t), d_n(t) \to 0$ as $n\to\infty$.
On the other hand, for any $T>0$  our assumptions imply (see \cite{DP1}) that \beq\label{AW}\sup_{n\ge 1}\sup_{0\leq t\leq T}\E |X^n_t|^2 +\sup_{0\leq t\leq T}\E |X_t|^2 <\infty.\end{equation}  So, the function $$g(t):=\limsup_{n\to \infty}\E |X_t^n-X_t|^2,\ \  t\ge 0 $$ is locally bounded. We will complete the proof of the theorem by showing $g(t)=0$. Taking $\limsup$ in (\ref{2.11}) we obtain
\beg{eqnarray}\label{2.12}
g(t)&\leq& C\int_0^t(K_{\si}(t-s)+K_b(t-s))g(s)\d s.
\end{eqnarray}
Given any $\beta>0$. Multiplying (\ref{2.12}) by $\e^{-\beta t}$ and integrating from $0$ to $T$ we get
\beq\label{2.13}\beg{split}
\int_0^T g(t)\e^{-\beta t}\d t&\leq C\int_0^T\d t \e^{-\beta t}\int_0^t(K_{\si}(t-s)+K_b(t-s))g(s)\d s \\
&=  C\int_0^T\e^{-\beta s}g(s)\d s\int_s^T(K_{\si}(t-s)+K_b(t-s))\e^{-\beta (t-s)}\d t \\
&= C\int_0^T\e^{-\beta s}g(s)\d s\int_0^{T-s}(K_{\si}(u)+K_b(u))\e^{-\beta u}\d u \\
&\leq \bigg(C\int_0^{T}(K_{\si}(u)+K_b(u))\e^{-\beta u}\d u\bigg)\int_0^T\e^{-\beta s}g(s)\d s.
\end{split}\end{equation}
Choosing $\beta>0$ sufficiently big so that $ C\int_0^{T}(K_{\si}(u)+K_b(u))\e^{-\beta u}\d u <1$, we deduce from (\ref{2.13}) that $\int_0^T g(t)\e^{-\beta t}\d t=0$ and  hence $g(t)=0,$ a.e.  By virtue of  (\ref{2.12}), we further conclude $g(t)=0$ for every $t\in [0,T]$, and thus finish the proof.
\end{proof}

\section{Proof of Theorem \ref{T1.1}}

According to Theorem \ref{TF} and using the monotone class theorem, it would be sufficient to prove Theorem 1.1  for the finite-dimensional setting, i.e. to prove the following result.

\beg{thm}\label{T2.1} Let $\H=\R^n$ and assume that {\bf (A1)} and {\bf (A3)} hold.
\beg{enumerate} \item[$(1)$] For any $f\in C_b^1(\R^n)$,   $$|\nn P_t f|^2\le 6^{1+\ff t {t_0}} P_t|\nn f|^2,\ \ t\ge 0.$$ \item[$(2)$] If  {\bf (A2)}  holds, then for any strictly positive $f\in \B_b(\R^n)$,
$$ P_t\log f(y)\le \log P_t f(x) +\ff{3\log 6}{\ll(\si)t_0( 1-6^{-\ff t{t_0}} )}|x-y|^2,\ \ x,y\in\R^n, t>0.$$
\item[$(3)$]  If  {\bf (A2)}  holds, then for any $f\in\B_b(\R^n)$,
$$|\nn P_tf|^2\le \ff{3\log 6}{t_0\ll(\si)(1- 6^{-\ff t{t_0}} )}\big\{P_t f^2-(P_tf)^2\big\},\ \ t>0.$$
\item[$(4)$] If          $|\si v|^2\le \bar\ll(\si)|v|^2$   holds for some constant $\bar\ll(\si)>0$ and all $v\in\R^n$, then
$$P_t f^2 -(P_t f)^2 \le \ff{12\bar\ll(\si) t_0(6^{\ff t{t_0}}-1)}{\log 6} P_t|\nn f|^2,\ \ f\in C_b^1(\R^n), t\ge 0.$$ \end{enumerate} \end{thm}

\beg{proof} In the present finite-dimensional setting,   {\bf (A1)} and {\bf (A3)} imply that $b$ and $\si$ are Lipschitz continuous. By a standard approximation argument we may and do assume that they are smooth with bounded gradients,  such that
\beq\label{A1} |T_s \nn_v b|^2\le K_b(s)|v|^2,\ \|T_s\nn_v \si\|^2\le K_\si(s)|v|^2,\ \ \ s>0, v\in\R^n, \end{equation}where $\nn_v$ denotes the directional derivative along $v$.
In this case the derivative process
$$\nn_v X_t:= \lim_{\vv\downarrow 0} \ff{X^{\cdot+\vv v}_t- X_t}\vv$$ solves the equation
$$\d \nn_v X_t= \big\{A\nn_v X_t+ (\nn_{\nn_vX_t}b)(X_t)\big\}\d t + (\nn_{\nn_v X_t}\si)(X_t)\d W_t,\ \ \nn_vX_0=v.$$
Since $\nn b$ and $\nn \si$ are bounded, this implies that
$$\sup_{s\in [0,t]} \E |\nn_vX_s|^2<\infty,\ \ t\ge 0.$$ We aim to find an upper bound of $\E |\nn_v X_t|^2$ independent of the dimension $n$ so that it can be passed to the infinite-dimensional setting.  To this end, let us observe that for any $s_0\ge 0$ we have
$$\nn_v X_t = T_{t-s_0} \nn_v X_{s_0} + \int_{s_0}^t T_{t-s} (\nn_{\nn_v X_s}b)(X_s) \d s + \int_{s_0}^t T_{t-s} (\nn_{\nn_v X_s}\si)(X_s)\d W_s,\ \ t\ge s_0.$$ Combining this with (\ref{A1}) we obtain
\beg{equation*}\beg{split} \E |\nn_v X_t|^2 &\le 3 \E|\nn_v X_{s_0}|^2 + 3  \int_{s_0}^tK_b(s-s_0) \E|\nn_v X_s|^2\d s + 3  \int_{s_0}^tK_\si(s-s_0) \E|\nn_v X_s|^2\d s\\
&\le   3 \E|\nn_v X_{s_0}|^2+ \big\{3\phi_b(t-s_0)+3\phi_\si (t-s_0)  \big\}\sup_{s\in[s_0,t]} \E|\nn_v X_s|^2\end{split}\end{equation*} for $t\ge s_0.$
Since the resulting  upper bound is increasing in $t\ge s_0$, it follows that
$$\sup_{s\in[s_0,t]} \E|\nn_v X_s|^2\le   3 \E|\nn_v X_{s_0}|^2+\big\{3\phi_b(t-s_0)+3\phi_\si (t-s_0)  \big\}\sup_{s\in[s_0,t]} \E|\nn_v X_s|^2 $$ holds for $t\ge s_0.$ Taking $t=s_0+t_0$ in this inequality leads to
$$\sup_{s\in [s_0,s_0+t_0]}\E|\nn_v X_s|^2\le 6 \E|\nn_v X_{s_0}|^2,\ \ s_0\ge 0.$$ Therefore,
\beq\label{GG} \E|\nn_v X_t|^2\le 6^{\ff{t+t_0}{t_0}} |v|^2,\ \ t\ge 0, v\in\R^n.\end{equation}
With this estimate in hand, we are able to complete the proof easily.

(1) follows from   (\ref{GG}) and the Schwarz inequality, more precisely
$$|\nn_v P_t f|^2 = \big|\nn_v \E f(X_t)\big|^2 =\big|\E\<\nn f(X_t), \nn_v X_t\>\big|^2 \le 6^{\ff{t+t_0}{t_0}} |v|^2P_t|\nn f|^2.$$

(2) follows from (1) according to the argument in \cite{RW10}, see Proposition \ref{P2.2} below for details.

(3) follows from   (1) by noting that
$$\ff{\d}{\d s}P_s(P_{t-s}f)^2 = 2P_s|\si\nn P_{t-s}f|^2 \ge 2\ll(\si) P_s|\nn P_{t-s}f|^2\ge 2\ll(\si) 6^{-\ff{s+t_0}{t_0}}|\nn P_t f|^2,\ \ s\in [0,t].$$

(4) follows from (1) since
$$\ff{\d}{\d s}P_s(P_{t-s}f)^2 = 2P_s|\si\nn P_{t-s}f|^2\le 2\bar\ll(\si)P_s|\nn P_{t-s}f|^2\le 2\bar\ll(\si) 6^{\ff{t-s-t_0}{t_0}}P_t|\nn f|^2,\ \ s\in [0,t].$$ \end{proof}

\beg{prp} \label{P2.2} If there exists a positive function $\Phi\in C([0,\infty))$ such that
\beq\label{GGE} |\nn P_t f|^2\le \Phi(t)P_t|\nn f|^2,\ \ t\ge 0, f\in C_b^1(\R^n),\end{equation} then for any strictly positive $f\in\B_b(\R^n)$,
$$P_t\log f(x)\le \log P_t f(y) +\ff{|x-y|^2}{2\ll(\si)\int_0^t \Phi(s)^{-1} \d s},\ \ t>0, x,y\in \R^n.$$\end{prp}

\beg{proof} For fixed $x,y\in\R^n, t>0$ and $h\in C^1([0,t];\R)$ with $h_0=0$ and $h_t=1$, let $x_s= (x-y)h_s+y,\ s\in [0,t].$  Combining (\ref{GGE}),
{\bf (A2}) and (see (2.3) in \cite{RW10})
$$\ff{\d }{\d s} P_s\log P_{t-s} f= -\ff 1 2 P_s |\si \nn \log P_{t-s} f|^2,$$ we obtain
\beg{equation*}\beg{split} & P_t\log f(x)-\log P_tf(y)= \int_0^t \ff{\d}{\d s } \left[\big(P_s\log P_{t-s}f\big)(x_s)\right]\d s\\
&= \int_0^t \Big\{h_s'\<x-y, \nn P_s\log P_{t-s}f\>- \ff 1 2 P_s |\si\nn\log P_{t-s}f|^2\Big\} (x_s)\d s\\
&\le \int_0^t \Big\{|h_s'|\cdot |x-y|\cdot |\nn P_s\log P_{t-s}f|- \ff {\ll(\si)}{ 2 \Phi(s)} |\nn P_s  \log P_{t-s}f|^2\Big\} (x_s)\d s\\
&\le \ff{ |x-y|^2}{2\ll(\si)} \int_0^t \Phi(s) (h_s')^2\d s.\end{split}\end{equation*} Taking
$$h_s= \ff{\int_0^s \Phi(u)^{-1}\d u}{\int_0^t \Phi(u)^{-1}\d u},\ \ s\in [0,t],$$ we complete the proof. \end{proof}

\beg{proof}[Proof of Theorem \ref{T1.1}] Let $P_t^n$ be the semigroup for $X_t^n$ solving the equation (\ref{2.2}). By Theorem \ref{TF} we have
\beq\label{*} P_tf(x)= \lim_{n\to\infty} P_t^nf(\scr P_nx),\ \ t\ge 0, f\in C_b(\H). \end{equation}

  Let $f\in C_b^1(\H)$. It is easy to see that if {\bf (Ai)} $(1\le {\bf i}\le 3)$ holds for $A,\si,b$ on $\H$, it also holds for $A_n,\si_n,b_n$ on $\H_n$. By Theorem \ref{T2.1} we have
$$\ff{|P_t^nf(\scr P_nx)-P_t^n f(\scr P_ny)|^2}{|x-y|^2} \le 6^{1+\ff t{t_0}} \int_0^1 (P_t^n |\nn f|^2)(s\scr P_n x+(1-s)\scr P_ny)\d s,\ \ x\ne y.$$ Letting $n\to\infty$ and using (\ref{*}), we arrive at
$$\ff{|P_tf(x)-P_t f(y)|^2}{|x-y|^2} \le 6^{1+\ff t{t_0}} \int_0^1 (P_t |\nn f|^2)(sx+(1-s)y)\d s,\ \ x\ne y.$$ This is equivalent to the gradient inequality in (1).
(4) can be proved similarly. By the same reason, it is easy to see that the inequalities in (2) and (3) hold for $f\in C_b(\H)$. Noting that  the inequality in (3) is equivalent to
$$ \ff{|P_tf(x)-P_t f(y)|^2}{|x-y|^2}\le \ff{3\log 6}{t_0\ll(\si) (1-6^{-\ff t{t_0}})}\int_0^1\{P_t f^2-(P_tf)^2\}(sx+(1-s)y)\d s,\ \ x\ne y,$$
by   the monotone class theorem, if  inequalities in (2) and (3) hold for all $f\in C_b(\H)$, they also hold for all $f\in \B_b(\H).$ \end{proof}

\section{Application to  white noise driven SPDEs}
In this section, we will apply our results to stochastic reaction-diffusion equations driven by space-time white noise
which are extensively studied in the literature, see \cite{DP2} and references therein.

  Consider the stochastic reaction-diffusion equation on a bounded closed domain $D\subset \R^d (d\ge 1)$:
\beq\label{4.1}\beg{cases}
&\frac{\partial u_t(\xi )}{\partial t} =  -(-\DD)^\aa u_t(\xi)
 +\psi(u_t(\xi))+ \phi(u_t(\xi
))\frac{\partial^{1+d}}{\partial t\partial \xi_1\cdots\pp \xi_d}W(t,\xi ),\\
 &u_0 =  g, \quad u_t|_{\pp D}=0, \quad \xi=(\xi_1,\cdots,\xi_d)\in D,\end{cases}
\end{equation}
where $\aa>0$ is a constant, $W(t,\xi )$ is a Brownian sheet on $\R^{d+1}$,  $\DD$ is the Dirichlet Laplacian on $D$, and $\phi,\psi$ are Lipschitz functions on $\R$, i.e. there exists a constant $C>0$ such that
\beq\label{LIP}
 |\psi(r)-\psi(s)|\leq c |r-s|, \quad |\phi(r)-\phi(s)|\leq c |r-s|,\ \ r,s\in \R.
 \end{equation}

The equivalent integral equation is (see \cite{W})
 \begin{eqnarray}
 u_t(\xi) &=& T_tg(\xi) +\int_{[0,t]\times D}T_{t-s}(\xi,\eta)\psi(u_s(\eta))\d s\d\eta  \nonumber\\
 &&+ \int_{[0,t]\times D}T_{t-s}(\xi,\eta)\phi(u_s(\eta))W(\d s, \d\eta),\ \ t\ge 0,
\end{eqnarray}
where $T_t$, $T_t(x,\eta)$ are the semigroup and  the heat kernel associated with the Dirichlet
Laplacian on $D$.

Top apply our main results to the present model, we reformulate the equation by using  the cylindrical Brownian motion on $\H:=  L^2(D).$ Let $A=-(-\DD)^\aa$.  Then  $-A$ has discrete spectrum with eigenvalues
$\{\ll_n\}_{n\ge 1}$   satisfying
\beq\label{EIG} \ff {n^{2\aa/d}}{C}\le \ll_n\le Cn^{2\aa/d},\ \ n\ge 1\end{equation} for some constant $C>1.$ Let $\{ e_n\}_{n\ge 1}$ be the corresponding unit eigenfunctions.
Since $e_m$ is independent of $\aa$, letting $\aa=1$ and using the classical Dirichlet heat kernel bound, we obtain
 \beq\label{UPP} \|e_m\|_\infty=\e \|T_{\ll_m^{-1}}\e_m\|_\infty \le \e \|T_{\ll_m^{-1}}\|_{L^2\to L^\infty}\le c_1 \ll_m^{d/4}\le c_2 \ss m,\ \ m\ge 1\end{equation} for some constants $c_1,c_2>0.$

Now,  define a sequence of independent Brownian motions by
$$\beta_n(t)=\int_{[0,t]\times D} e_n(\eta)W(\d s, \d\eta), \quad  n\geq 1.$$
Then
 $$W_t:=\sum_{n=1}^{\infty} \beta_n(t)e_n $$ is a  cylindrical Brownian motion on $\H$. Let
$$b (u)(\xi)=\psi(u(\xi )),\ \ \{\sigma(u) x\} (\xi)=\phi(u(\xi))\cdot x(\xi),\ \ u, x\in \H, \xi\in D.$$ It is
easy to see that the reaction-diffusion diffusion equation (\ref{4.1}) can be reformulated as
$$\d u_t= Au_t\d t+b(u_t)\d t +\sigma(u_t)\d W_t.$$
Obviously, $\si$ takes values in the space of bounded linear operators on $\H$ if and only if $\phi$ is bounded. So, in general $\si$ is not a  Hilbert-Schmidt perturbation of  any bounded linear operator as indicated in the Introduction.

\beg{thm}\label{T4.1} Let $A, b$ and $\si$ be given above such that  $ \phi^2\ge \ll$ for some constant $\ll>0$. \beg{enumerate} \item[$(1)$] If $\aa>d$   then conditions {\bf (A1)}-{\bf (A4)} hold for $\ll(\si)=\ll, K_b\equiv C$ and
\beq\label{KA}  K_\si(t)=C \sum_{m=1}^\infty m\e^{-\dd t m^{2\aa/d}}\end{equation}
for some constants $C,\dd>0$,
so that Theorem $\ref{T1.1}$   and Corollary  $\ref{C1.2}$ apply to the semigroup associated to solutions of
equation $(\ref{4.1}). $
\item[$(2)$] Let $D=\prod_{i=1}^d [a_i, b_i]$ for some $b_i>a_i, 1\le i\le d$. If $\aa>\ff d 2$ then {\bf (A1)}-{\bf (A4)} hold for $\ll(\si)=\ll, K_b\equiv C $ and
\beq\label{KB} K_\si(t)=C \sum_{m=1}^\infty  \e^{-\dd t m^{2\aa/d}}\end{equation} for some constants $C,\dd>0, $ so that Theorem $\ref{T1.1}$   and Corollary  $\ref{C1.2}$ apply to the semigroup associated to solutions of
equation $(\ref{4.1}). $
\item[$(3)$] In the situations of $(1)$ and $(2)$, there  exists a constant $\vv_0>0$ such that $P_t$ has a unique invariant probability measure provided
\beq\label{**DW} |\phi(s)|^2+|\psi(s)|^2\le \vv_0 |s|^2 +C_0,\ \ s\in\R\end{equation} holds for some constant $C_0>0.$ \end{enumerate}
\end{thm}

\beg{proof} Since {\bf (A4)} is obvious due to (\ref{EIG}) and {\bf (A2)} with $\ll(\si)=\ll$ follows from $\phi^2\ge\ll$, it suffices to
verify {\bf (A1)} and {\bf (A3)} for the desired $K_b$ and $K_\si$.
By the contraction of $T_t$ and (\ref{LIP}), we have
$$|T_t(b(x)-b(y))|\leq c|x-y|, \quad x, y\in \H,$$ where $|\cdot|$ is now the $L^2$-norm on $D$.
Then {\bf (A1)} holds for   $K_b\equiv c^2$.

Below, we verify {\bf (A3)} and the existence of the invariant probability measure respectively.

(1) By the definition of $\si,$ (\ref{LIP}) and (\ref{UPP}), we have
\beq\label{PPQ} \beg{split} &\|T_t(\sigma(x)-\sigma(y))\|_{HS}^2=\sum_{n=1}^{\infty}|T_t(\sigma(x)-\sigma(y)) e_n |^2\\
&=\sum_{n=1}^{\infty}\sum_{m=1}^{\infty}\<T_t(\sigma(x)-\sigma(y)) e_n, e_m\>^2 \\
 & =\sum_{m=1}^{\infty}\e^{-2t\ll_m }\sum_{n=1}^{\infty}\<(\sigma(x)-\sigma(y)) e_n,
e_m\>^2\\
&=\sum_{m=1}^{\infty}\e^{-2t\ll_m }|(\sigma(x)-\sigma(y))^* e_m |^2 \\
& =\sum_{m=1}^{\infty}\e^{-2t\ll_m }\int_D
\big|(\phi(x(\xi))-\phi(y(\xi))) e_m(\xi)\big|^2 \d\xi\\
& \leq c^2  |x-y |^2 \sum_{m=1}^\infty \|e_m\|_\infty^2    \e^{-\dd t m^{2\aa/d}}.
\end{split}\end{equation}  for some constant $\dd>0$. Combining this with (\ref{UPP}) we obtain
$$\|T_t(\sigma(x)-\sigma(y))\|_{HS}^2\le
  C\sum_{m=1}^\infty m\,\e^{-\dd t m^{2\aa/d}}$$ for some constant $C>0$. Moreover,
\beq\label{LM} \int_0^t \|T_s \si(0)\|_{HS}^2=  \phi(0)^2 \int_0^t \|T_s\|_{HS}^2 \d s\\
 \le   C' \sum_{m=1}^\infty \e^{-\dd t m^{2\aa/d}} \end{equation} holds for some constant $C'>0$. Therefore, if $\aa>d$ then {\bf(A3)} holds for $K_\si$ given in (\ref{KA}) since in this case
$$\int_0^\infty \sum_{m=1}^\infty m\e^{-\dd t m^{2\aa/d}}\d t =\ff 1 \dd \sum_{m=1}^\infty \ff 1 {m^{(2\aa-d)/d}}<\infty.$$

(2) When $D=\prod_{i=1}^d [a_i, b_i]$ for some $b_i>a_i, 1\le i\le d$,  the eigenfunctions $\{e_m\}_{m\ge 1}$ are uniformly bounded, i.e. $\|e_m\|_\infty\le C$ holds for some constant $C>0$ and all $m\ge 1$. Combining this with (\ref{PPQ}),  we obtain
$$\|T_t(\sigma(x)-\sigma(y))\|_{HS}^2\le C|x-y|^2 \, \sum_{m=1}^\infty  \e^{-\dd t m^{2\aa/d}}$$ for some constants $C,\dd>0$. Combining this with (\ref{LM}), we conclude that
{\bf(A3)} holds for $K_\si$ given in (\ref{KB}) provided $\aa>\ff d 2.$

(3) The uniqueness of the invariant probability measure follows from Corollary \ref{C1.2}(3), it suffices to prove the existence by verifying conditions (i)-(iv) in \cite[Theorem 6.1.2]{DP2}.  By (\ref{LIP}), (\ref{EIG}) and {\bf (A3)}, conditions (i) and (iii) hold. It remains to verify condition (ii), i.e.
\beq\label{**DD} \int_0^1 s^{-\vv} K_\si(s)\d s<\infty \ \text{for\ some\ }\vv\in (0,1);\end{equation} and condition (iv), which is implied by
\beq\label{**DD2} \sup_{t\ge 0} \E|u_t|^2<\infty.\end{equation}

Let $K_\si$ be in (\ref{KA}) with $\aa>d$. Then for $\vv\in (0, \ff{\aa-d}\aa),$
\beg{equation*}\beg{split} &\int_0^1 s^{-\vv} K_\si(s)\d s = C\sum_{m=1}^\infty m \int_0^1 s^{-\vv} \e^{-\dd s m^{2\aa/d}}\d s\\
&\le C\sum_{m=1}^\infty m\bigg(\int_0^{m^{-2\aa/d}} s^{-\vv}\d s + m^{2\aa\vv/d} \int_{m^{-2\aa/d}}^1 \e^{-\dd s m^{2\aa/d}}\d s\bigg)\\
&\le C(\vv) \sum_{m=1}^\infty m^{1-2(1-\vv)\aa/d}<\infty,\end{split}\end{equation*} where $C(\vv)>0$ is a constant depending on $\vv$. Similarly, (\ref{**DD}) holds for $K_\si$ in (\ref{KB}) with $\aa>\ff d 2$ and $\vv\in (0, \ff{2\aa-d}{2\aa}).$

Next, by (\ref{**DW}) we have
\beg{equation*}\beg{split} \E|u_t|^2 &\le C_1 |g|^2 +C_1 \int_0^t \sum_{m=1}^\infty \|e_m\|_\infty^2 e^{-2\lambda_m(t-s)}\big(C_0+\vv_0 \E|u_s|^2\big)\d s\\
&\le C_1|g|^2 +C_2 \Big(C_0+\vv_0\sup_{s\in [0,t]}\E|u_s|^2\Big) \int_0^tK_\si(s)\d s,\ \ t\ge 0\end{split}\end{equation*} for some constants $C_1,C_2>0,$ where $K_\si$ is in (\ref{KA}) with $\aa>d$ or in (\ref{KB}) with $\aa>\ff \aa 2$ such that $\int_0^\infty K_\si(s)\d s <\infty$ as observed above. So, there exist   constants  $C_3,C_4>0$ such that
$$\E |u_t|^2\le C_3 +C_4\vv_0\sup_{s\in [0,t]}\E |u_s|^2,\ \ t\ge 0.$$ Taking   $\vv_0=\ff 1{2C_4},$ we obtain
$$\sup_{s\in [0,t]}\E|u_s|^2 \le C_3+ C_4\vv_0 \sup_{s\in [0,t]} \sup_{r\in [0,s]} \E|u_r|^2= C_3 +\ff 1 2\sup_{s\in [0,t]}\E|u_s|^2,\ \ t\ge 0.$$  Since by (\ref{AW}) $\sup_{s\in [0,t]}\E|u_s|^2<\infty$, this implies (\ref{**DD2}).\end{proof}

\beg{thebibliography}{99}
 \bibitem{ATW06} M. Arnaudon, A. Thalmaier, F.-Y. Wang,
  \emph{Harnack inequality and heat kernel estimates
  on manifolds with curvature unbounded below,} Bull. Sci. Math.   130(2006), 223--233.

\bibitem{ATW09} M. Arnaudon, A. Thalmaier, F.-Y. Wang,
  \emph{Gradient estimates and Harnack inequalities on non-compact Riemannian manifolds,}
   Stoch. Proc. Appl.   119(2009), 3653--3670.

   \bibitem{ATW09} M. Arnaudon, A. Thalmaier, F.-Y. Wang,
  \emph{Gradient estimates and Harnack inequalities on non-compact Riemannian manifolds,}
   Stoch. Proc. Appl. 119(2009), 3653--3670.

\bibitem{DP1} G.D.~Prato, J.~Zabczyk, \emph{
Stochastic Equations in Infinite Dimensions,}
Cambridge University Press, 1992.

\bibitem{DP2} G.D.~Prato, J.~Zabczyk, \emph{
Ergodicity for  Infinite Systems,}
Cambridge University Press, 1996.

\bibitem{GW} A. Guillin, F.-Y. Wang,
  \emph{Degenerate Fokker-Planck equations : Bismut formula, gradient estimate  and Harnack inequality,}  J. Diff. Equat. 253(2012), 20--40.

\bibitem{DP} C. Donati-Martin, J. Zabczyk, \emph{White noise driven SPDEs with reflection,} Probab. Theory Relat. Fields 95(1993), 1--24.

\bibitem{L} W. Liu, \emph{Harnack inequality and applications for stochastic evolution equations with monotone drifts,}  J. Evol. Equ. 9 (2009),  747--770.

\bibitem{LW} W.  Liu, F.-Y.  Wang,  \emph{Harnack inequality and strong Feller property for stochastic fast-diffusion equations,} J. Math. Anal. Appl. 342(2008), 651--662.

\bibitem{O} S.-X. Ouyang, \emph{Harnack inequalities and applications for multivalued stochastic evolution equations,}  Inf. Dimen.
Anal. Quant. Probab. Relat. Topics.

\bibitem{PZ} S. Peszat, J. Zabaczyk, \emph{Strong Feller property and irreducibility for diffusions on Hilbert spaces,} Ann. Probab. 23(1995), 157--172.

\bibitem{RW10} M. R\"ockner, F.-Y. Wang, \emph{Log-Harnack  inequality for stochastic differential equations in Hilbert spaces and its consequences, } Infin. Dimens. Anal. Quant. Probab.  Relat. Topics 13(2010), 27--37.

\bibitem{W97} F.-Y. Wang,  \emph{Logarithmic Sobolev inequalities on noncompact Riemannian manifolds,} Probab. Theory Relat. Fields 109(1997), 417-424.
\bibitem{W10} F.-Y. Wang \emph{Harnack inequalities on manifolds with boundary and applications,}    J.
Math. Pures Appl.   94(2010), 304--321.

\bibitem{W11} F.-Y. Wang,  \emph{Harnack inequality for SDE with multiplicative noise and extension to
               Neumann semigroup on nonconvex manifolds,} Ann. Probab. 39(2011), 1447-1467.

\bibitem{WWX10}  F.-Y. Wang, J.-L. Wu , L. Xu, \emph{Log-Harnack inequality for stochastic Burgers equations and
applications},     J. Math. Anal. Appl. 384(2011), 151--159.

\bibitem{WY11} F.-Y. Wang, C. Yuan, \emph{Harnack inequalities for functional SDEs with multiplicative noise and applications,}   Stoch. Proc. Appl.    121(2011), 2692--2710.

 \bibitem{WZ} F.-Y. Wang, T. Zhang, \emph{Gradient estimates for stochastic evolution equations with non-Lipschitz coefficients,}  J. Math. Anal. Appl. 365(2010), 1--11.

 \bibitem{W} J. Wash, \emph{ An introduction to stochastic partial differential equations,} Lecture Notes in Mathematics 1180 (1984) 265--439.

\bibitem{X} L. Xu, \emph{A modified log-Harnack inequality and asymptotically strong Feller property,}  J. Evol. Equ. 11(2011),  925--942.

\bibitem{ZS} S.-Q.  Zhang, \emph{Harnack inequality for semilinear SPDE with
multiplicative noise,} arXiv:1208.2343.

\bibitem{ZT} T. Zhang, \emph{White noise driven SPDEs with reflection: strong Feller properties and Harnack inequalities,} Potential Anal. 33(2010), 137--151.

\end{thebibliography}
\end{document}